\documentclass[12pt]{article}
\usepackage[utf8]{inputenc}
\usepackage{amssymb,amsmath,amsfonts,eucal,mathrsfs,amsthm} 
\usepackage[colorinlistoftodos]{todonotes}
\usepackage[allcolors=blue]{hyperref}

\usepackage{empheq}
\newtheorem{theorem}{Theorem}
\newtheorem{proposition}[theorem]{Proposition}
\newtheorem{lemma}[theorem]{Lemma}

\newtheorem{corollary}[theorem]{Corollary}

\newtheorem{remark}[theorem]{Remark}

\newtheorem{examples}[theorem]{Examples}

\theoremstyle{definition}
\newcommand{\R}{\mathbb{R}}

\newcommand{\Sf}{\mathbb{S}}

\newcommand{\spa}{\mbox{span}}
\newcommand{\hess}{\mbox{Hess\,}}

\newcommand{\nap}{\nabla^{\perp}}
\newcommand{\nab}{\tilde\nabla}

\newcommand{\End}{\mbox{End}}
\newcommand{\Hom}{\mbox{Hom}}

\newcommand{\Y}{\mathcal{Y}\,}
\newcommand{\D}{\mathcal{D}}
\newcommand{\Les}{\mathbb{L}}
\def\<{{\langle}}
\def\>{{\rangle}}

\def\Sal{{\cal S}}

\def\T{{\cal T}}

\def\Y{{\cal Y}}
\def\n{\nabla}
\def\d{\partial}
\def\a{\alpha}

\def\be{\begin{equation} }
\def\ee{\end{equation} }

\def\proof{\noindent{\it Proof:  }}
\def\qed{\ifhmode\unskip\nobreak\fi\ifmmode\ifinner
\else\hskip5 pt \fi\fi\hbox{\hskip5 pt \vrule width4 pt
height6 pt  depth1.5 pt \hskip 1pt }}
\makeatletter
\newcommand{\subjclass}[2][]{\let\@oldtitle\@title
\gdef\@title{\@oldtitle\footnotetext{#1 
\emph{Mathematics Subject Classification:} #2}}}
\newcommand{\keywords}[1]{\let\@@oldtitle\@title
\gdef\@title{\@@oldtitle\footnotetext
{\emph{Key words and phrases.} #1.}}}
\makeatother
\begin{document}

\title{Conformal infinitesimal variations\\ of 
submanifolds}
\author{M. Dajczer and M. I. Jimenez}
\date{}
\maketitle

\begin{abstract} This paper belongs to the realm of conformal 
geometry and deals with Euclidean submanifolds that admit  
smooth variations that are infinitesimally conformal. Conformal 
variations of Euclidean submanifolds is a classical subject 
in differential geometry. In fact, already in 1917 Cartan
classified parametrically the Euclidean hypersurfaces that 
admit nontrivial conformal variations. Our first main result 
is a Fundamental theorem for conformal 
infinitesimal variations. The second is a rigidity theorem 
for Euclidean submanifolds that lie in low codimension.
\end{abstract}

This paper deals with a subject in conformal geometry, namely, 
with smooth conformal infinitesimal variations of Euclidean 
submanifolds of any dimension and in any codimension. 
Until now the study of this class of variations has received 
limited attention; see \cite{Ya} for an exception.  This is 
certainly not the situation for the more restricted case of 
isometric infinitesimal variations. In fact, for hypersurfaces 
the study of these variations  is a classical subject already 
considered by Sbrana \cite{Sb} at the beginning of the $20^{th}$ 
century after the earlier rigidity result contained in Cesàro's 
book \cite{Ce} from 1896. For recent results on the subject, we 
refer to \cite{DV} and \cite{Ji} in the hypersurface case and to 
\cite{DJ} and \cite{DJ2} for submanifolds in higher codimension.

First Sbrana \cite{Sb2} and subsequently Cartan \cite{Ca0} classified 
parametrically Euclidean hypersurfaces that admit nontrivial isometric 
variations. 
Shortly after, Cartan in \cite{Ca} classified the Euclidean hypersurfaces 
that admit nontrivial conformal variations; see also \cite{DT0} 
and \cite{DT}. A conformal variation of a given isometric immersion 
$f\colon M^n\to\R^m$ of a Riemannian manifold $(M^n, \<\,,\,\>)$ 
into Euclidean space is a smooth variation $F\colon I\times M^n\to\R^m$, 
where $0\in I\subset\R$ is an open interval and $f_t=F(t,\cdot)$ with 
$f_0=f$ is a conformal immersion for any $t\in I$. Hence, there 
is a positive function $\gamma\in C^\infty(I\times M^n)$ with 
$\gamma(0,x)=1$ such that 
\be\label{gamma}
\gamma(t,x)\<f_{t*}X,f_{t*}Y\>=\<X,Y\>
\ee
for any tangent vector fields $X,Y\in\mathfrak{X}(M)$. Here and 
in the sequel, we use the same notation for the inner products in 
$\R^m$ and $M^n$, and denote by $\nab$ and $\n$ the respective 
associated Levi-Civita connections. The derivative of \eqref{gamma} 
computed at $t=0$ gives that the variational vector field 
$\T=F_*\d/\d t|_{t=0}$ of $F$ has to satisfy the condition
\be\label{cib0}
\<\nab_X\T,f_*Y\>+\<f_*X,\nab_Y\T\>=2\rho\<X,Y\>
\ee
where $\rho(x)=-(1/2)\d\gamma/\d t(0,x)$.

Trivial conformal variations are the ones induced by a 
composition of the immersion with a smooth family of conformal 
transformations of the Euclidean ambient space. 
In this case, the variational vector field is, at least locally, 
the restriction of a conformal Killing vector field of the ambient 
space to the submanifold. Recall that conformal transformations 
of Euclidean space are characterized by Liouville's classical 
theorem; see \cite{To} for a nice discussion of this result.
\vspace{1ex}

In this paper, we deal with the weaker concept
of conformal infinitesimal  variation, that is,
the infinitesimal analogue of a conformal variation.
\vspace{1ex}

A smooth variation $F\colon I\times M^n\to\R^m$ of a given isometric
immersion $f\colon M^n\to\R^m$ is called a 
\emph{conformal infinitesimal variation} 
if there is a function $\gamma\in C^\infty(I\times M^n)$ satisfying 
$\gamma(0,x)=1$ such that
\be\label{varcond}
\frac{\d}{\d t}|_{t=0}\left(\gamma(t,x)\<f_{t*}X,f_{t*}Y\>\right)=0
\ee
for any $X,Y\in\mathfrak{X}(M)$. 
\vspace{1ex}

It is already well-known from classical differential geometry 
of submanifolds that the appropriate approach to study infinitesimal 
variations is to deal  with the variational vector field. In our 
case, that this is the way to proceed is justified in the sequel.
\vspace{1ex}

We have from \eqref{varcond} that the variational vector field $\T$ 
of $F$ has to satisfy condition \eqref{cib0}. This leads to the 
following definition.
\vspace{1ex}

A \emph{conformal infinitesimal bending} with \emph{conformal 
factor} $\rho\in C^\infty(M)$ of an isometric immersion 
$f\colon M^n\to\R^m$ of a Riemannian manifold $M^n$ into 
Euclidean space is a smooth section $\T\in\Gamma(f^*T\R^m)$ 
that satisfies 
\be\label{cib}
\<\nab_X\T,f_*Y\>+\<f_*X,\nab_Y\T\>=2\rho\<X,Y\>
\ee
for any  $X,Y\in\mathfrak{X}(M)$.
\vspace{1ex}

On one hand, there is a conformal infinitesimal 
bending associated to any conformal infinitesimal variation.
On the other hand, associated to a conformal infinitesimal 
bending we have the variation 
$F\colon\R\times M^n\to\R^m$ given by
\be\label{unique}
F(t,x)=f(x)+t\T(x).
\ee
This is a conformal infinitesimal  variation with variational
vector field $\T$ since \eqref{varcond} is satisfied for 
$\gamma(t,x)=e^{-2t\rho(x)}$. By no means
\eqref{unique} is unique with this property, although it may 
be seen as the simplest one.  In fact, new conformal 
infinitesimal variations with variational vector field $\T$ 
are obtained by adding to \eqref{unique} terms of the type 
$t^k\delta$, $k>1$, where $\delta\in\Gamma(f^*T\R^{n+1})$ and, 
maybe, for restricted values of the parameter $t$.

We observe that the notion of conformal infinitesimal variation 
is indeed a concept in conformal geometry. In fact, let 
$F\colon I\times M^n\to\R^m$ be a conformal infinitesimal variation
of $f\colon M^n\to\R^m$.  Then, let $G\colon I\times M^n\to\R^m$ be 
the variation given by $G=\psi\circ F$ where $\psi$ is a conformal 
transformation of $\R^m$ with positive conformal factor 
$\lambda\in C^\infty(\R^m)$. We claim that $G$ is a conformal 
infinitesimal variation of $g=\psi\circ f$ where
$$
\tilde{\gamma}(t,x)=\gamma(t,x)-2t\<\T(x),\nab\log\lambda(f(x))\>.
$$
In fact, we have using \eqref{varcond} that 
\begin{align*}
&\frac{\d}{\d t}|_{t=0}\left(\tilde{\gamma}\< g_{t*}X,g_{t*}Y\>\right)\\
&=\frac{\d}{\d t}|_{t=0}\left((\gamma(t,x)-2t\<\T(x),\nab\log\lambda(f(x))\>)
\lambda^2(F(t,x))\< f_{t*}X,f_{t*}Y\>\right)\\
&=\< f_*X,f_*Y\>\frac{\d}{\d t}|_{t=0}(\lambda^2(F(t,x)))
-2\frac{\d}{\d t}|_{t=0}(t\lambda^2\<\T,\nab\log\lambda\>
\< f_{t*}X,f_{t*}Y\>)\\
&=2\lambda\<\T,\nab\lambda\>\< f_*X,f_*Y\>
-2\lambda^2\<\T,\nab\log\lambda\>\< f_*X,f_*Y\>=0,
\end{align*}
and this proves the claim.
\vspace{1ex}

In view of the above, we say that a conformal infinitesimal 
variation of $f\colon M^n\to\R^m$ is a \emph{trivial conformal 
infinitesimal variation} if the associated conformal infinitesimal 
bending is trivial. In turn, that a conformal infinitesimal bending
is \emph{trivial} means that it is locally the restriction  of a 
conformal Killing vector field of the Euclidean ambient space
to the submanifold. Finally, if any conformal infinitesimal bending 
of $f$ is trivial we say that the submanifold is \emph{conformally 
infinitesimally rigid}.
\vspace{1ex}

For a conformal infinitesimal bending $\T\in\Gamma(f^*T\R^m)$ 
with conformal factor $\rho\in C^\infty(M)$
of an isometric immersion $f\colon M^n\to\R^m$, we  first 
show that $\T$ together with the second fundamental form 
$\a\colon TM\times TM\to N_fM$ of $f$ determine an 
\emph{associate pair} of tensors $(\beta,{\cal E})$,
where $\beta\colon TM\times TM\to N_fM$ is symmetric and 
${\cal E}\colon TM\times N_fM\to N_fM$ satisfies the 
compatibility condition
\be\label{anti}
\<{\cal E}(X,\eta),\xi\>+\<{\cal E}(X,\xi),\eta\>=0
\ee
for any $X\in\mathfrak{X}(M)$ and $\eta,\xi\in\Gamma(N_fM)$.
Subsequently, we prove that the pair $(\beta,{\cal E})$ satisfies 
the following fundamental system  of equations, where 
by the term fundamental we  mean that they are the integrability 
condition for the existence of a conformal infinitesimal bending. 
\begin{empheq}[left=\hypertarget{S}{(S)}\empheqlbrace]{align}
A_{\beta(Y,Z)}X+B_{\alpha(Y,Z)}X-A_{\beta(X,Z)}Y-B_{\alpha(X,Z)}Y\label{derGaussC}\nonumber\\
+(X\wedge HY-Y\wedge HX)Z=0\\
(\nabla^{\perp}_X\beta)(Y,Z)-(\nabla_Y^{\perp}\beta)(X,Z)
={\cal E}(Y,\a(X,Z))-{\cal E}(X,\a(Y,Z))\label{derCodazziC}\nonumber\\
+\,\<Y,Z\>\a(X,\nabla\rho)-\<X,Z\>\a(Y,\nabla\rho)\\
(\nap_X{\cal E})(Y,\eta)-(\nap_Y{\cal E})(X,\eta)
=\beta(X,A_\eta Y)-\beta(A_\eta X,Y)\label{derRicciC}\nonumber\\
+\,\a(X,B_\eta Y)-\a(B_\eta X,Y)
\end{empheq}
where $X,Y,Z\in\mathfrak{X}(M)$, $\eta\in\Gamma(N_fM)$ and 
$A_\eta,B_\eta, H\in\Gamma(\End(TM))$ are given by
$\<A_\eta X,Y\>=\<\a(X,Y),\eta\>$, 
$\<B_\eta X,Y\>=\<\beta(X,Y),\eta\>$ and $HX=\nabla_X\nabla\rho$ respectively.
\vspace{1ex}

An infinitesimal bending is a conformal infinitesimal 
bending with conformal factor $\rho=0$. It is said to be trivial 
if it is locally the restriction to the submanifold of a Killing 
vector field of the Euclidean ambient space. Let $\T_1$ be a 
conformal infinitesimal bending of  $f$ with 
conformal factor $\rho$ and let $\T_0$ be an infinitesimal 
bending of $f$. Then $\T_2=\T_1+\T_0$ satisfies \eqref{cib}, and 
therefore it is also a conformal infinitesimal bending of $f$ with 
conformal factor $\rho$.

The sum of any two conformal infinitesimal bendings is 
again a conformal infinitesimal bending. From now on, 
we identify two conformal infinitesimal bendings if they 
differ by a trivial conformal infinitesimal bending. 

\begin{theorem}\label{fundamental}
Let $f\colon M^n\to\R^m$, $n\geq 3$, be an isometric 
immersion of a simply connected Riemannian manifold. 
A triple $(\beta,{\cal E},\rho)\neq 0$ formed by a 
symmetric tensor $\beta\colon TM\times TM\to N_fM$, 
a tensor ${\cal E}\colon TM\times N_fM\to N_fM$ for which \eqref{anti} holds and $\rho\in C^\infty(M)$ that  
satisfies system $\hyperlink{S}{(S)}$ determines a 
unique conformal infinitesimal bending of $f$.
\end{theorem}

The above result takes a rather simpler form in the 
hypersurface case. In fact, let $f\colon M^n\to\R^{n+1}$ 
be a hypersurface with shape operator $A$ corresponding to the 
Gauss map $N\in\Gamma(N_fM)$. Associated to a conformal 
infinitesimal bending  we are now reduced to consider 
the tensor ${\cal B}\in\Gamma(\End(TM))$ given by 
$\beta(X,Y)=\<{\cal B}X,Y\>N$. Then the fundamental 
system of equations takes the form
\be\label{derGaussChip}
{\cal B}X\wedge AY-{\cal B}Y\wedge AX
+X\wedge HY-Y\wedge HX=0
\ee
and
\be\label{derCodazziChip}
(\nabla_X{\cal B})Y-(\nabla_Y{\cal B})X+(X\wedge Y)A\nabla\rho=0
\ee
for any $X,Y\in\mathfrak{X}(M)$. 

\begin{corollary}\label{corollary}
Let $f\colon M^n\to\R^{n+1}$, $n\geq 3$, be an isometric 
immersion of a simply connected Riemannian manifold. 
Then a symmetric tensor $0\neq{\cal B}\in\Gamma(\End(TM))$ 
and $\rho\in C^\infty(M)$ that satisfy \eqref{derGaussChip} 
and \eqref{derCodazziChip} determine a unique 
conformal infinitesimal bending of $f$.
\end{corollary}

The second main result in this paper is a rigidity theorem 
for conformal infinitesimal variations of submanifolds of low 
codimension.  The limitation on the codimension is due to 
the use of a result in the theory of flat bilinear forms
that is known to be false for higher codimensions.
\vspace{1ex}

The \emph{conformal $s$-nullity} $\nu_s^c(x)$  at $x\in M^n$, 
$1\leq s\leq m-n$, of an isometric immersion 
$f\colon M^n\to\R^m$ is defined as
$$
\nu_s^c(x)=\max\{\dim\mathcal{N}(\a_{U^s}-\<\,,\,\>\zeta)(x)
\colon U^s\subset N_fM(x)\;\mbox{and}\;\zeta\in U^s\}
$$
where $\a_{U^s}=\pi_{U^s}\circ\a$, $\pi_{U^s}\colon N_fM\to U^s$
is the orthogonal projection onto the $s$-dimensional subspace 
$U^s\subset N_fM$ and
$$
\mathcal{N}(\a_{U^s}-\<\,,\,\>\zeta)(x)=\{Y\in T_xM:
\a_{U^s}(Y,X)-\<Y,X\>\zeta=0\;\mbox{for all}\;X\in T_xM\}.
$$
The conformal $s$-nullity is a concept in conformal geometry 
since it is easily seen to be invariant under a conformal change 
of the metric of the ambient space.  
\vspace{1ex}

With respect to the next result, we observe that it has been 
shown in \cite{DV} that the set of Euclidean hypersurfaces 
admitting a nontrivial infinitesimal variation is 
much larger than the ones allowing an isometric variation. 
In \cite{DJV} it is shown that for hypersurfaces 
in the conformal case the situation is similar.

\begin{theorem}\label{rigidity}
Let $f\colon M^n\to \R^{n+p}$, $n\geq 2p+3$, be an isometric 
immersion with codimension $1\leq p\leq 4$. Assume that the 
conformal $s$-nullities of $f$ at any point of $M^n$ satisfy 
$\nu_s^c\leq n-2s-1$ for all $1\leq s\leq p$. 
Then $f$ is conformally infinitesimally rigid.
\end{theorem}

By the above result, in the case of hypersurfaces 
$f\colon M^n\to\R^{n+1}$, $n\geq 5$, the existence of a nontrivial 
conformal infinitesimal variation requires the presence at any point 
of a principal curvature of multiplicity at least $n-2$.

Theorem \ref{rigidity} is the version for conformal infinitesimal
variations of the rigidity result for conformal immersions due 
to do Carmo and Dajczer \cite{CD}, where the concept of conformal 
$s$-nullity has been introduced.  Moreover, the corresponding result 
for variations that are infinitesimally isometric was given
by Dajczer and Rodríguez \cite{DR}. Both results, as well as additional 
information, can be found in \cite{DT}.  With respect to the latter 
result, in sharp contrast with the situation in this paper a very short 
proof was possible in \cite{DT} by the use of a classical trick that 
fails completely in the conformal case.

\section{The fundamental equations}

In this section, we define a pair of tensors $(\beta,{\cal E})$ 
associated to a conformal infinitesimal bending $\mathcal{T}$
of $f\colon M^n\to\R^m$ and show that they satisfy the system 
of equations $\hyperlink{S}{(S)}$.  
\vspace{2ex}

Let $L\in\Gamma(\Hom(TM,f^*T\R^m))$ be the tensor defined by
$$
LX=\nab_X \T-\rho f_*X=\T_*X-\rho f_*X
$$
for any $X\in\mathfrak{X}(M)$. Then \eqref{cib} in terms 
of $L$ has the form
\be\label{inf}
\<LX,f_*Y\>+\<f_*X,LY\>=0
\ee
for all $X,Y\in\mathfrak{X}(M)$. Let 
$B\colon TM\times TM\to f^*T\R^m$ be the tensor given by
$$
B(X,Y)=(\nab_XL)Y=\nab_XLY-L\n_XY
$$
for any $X,Y\in\mathfrak{X}(M)$. Then the 
tensor
$\beta\colon TM\times TM\to N_fM$ is defined by
$$
\beta(X,Y)=(B(X,Y))_{N_fM}
$$
for any $X,Y\in\mathfrak{X}(M)$.  Flatness of the ambient space
and that
$$
\beta(X,Y)=(\nab_X\nab_Y\T-\nab_{\nabla_XY}\T)_{N_fM}-\rho\a(X,Y)
$$
give that $\beta$ is symmetric.

Let $\Y\in\Gamma(\Hom(N_fM, TM))$ be defined by
\be\label{Y}
\<\Y\eta,X\>+\<\eta,LX\>=0.
\ee
Then, let ${\cal E}\colon TM\times N_fM\to N_fM$ be the tensor 
given by
$$
{\cal E}(X,\eta)=\a(X,\Y\eta)+(LA_\eta X)_{N_fM}.
$$
We have
\begin{align*}
\<{\cal E}(X,\eta),\xi\>
&=\<\a(X,\Y\eta)+LA_\eta X,\xi\>
=\<A_\xi X,\Y\eta\>-\<\Y\xi,A_\eta X\>\\
&=-\<LA_\xi X,\eta\>-\<\a(X,\Y\xi),\eta\>
=-\<{\cal E}(X,\xi),\eta\>,
\end{align*}
and hence condition \eqref{anti} is satisfied.

\begin{lemma}\label{tam} We have that
$$
(B(X,Y))_{f_*TM}=f_*(\Y\a(X,Y)+(X\wedge\nabla\rho)Y)
$$
where $\nabla\rho$ denotes the gradient of $\rho$. 
\end{lemma}

\proof We have to show that
$$
C(X,Y,Z)=\<(B-f_*\Y\a)(X,Y),f_*Z\>+\<X,Y\>\<Z,\nabla\rho\>
-\<Y,\nabla\rho\>\<X,Z\>
$$
vanishes for any $X,Y,Z\in\mathfrak{X}(M)$. The 
derivative of \eqref{inf} gives
\begin{align*}
0&=\<\nab_ZLX,f_*Y\>+\<LX,\nab_Zf_*Y\>
+\<\nab_ZLY,f_*X\>+\<LY,\nab_Zf_*X\>\\ 
&=\<B(Z,X),f_*Y\>+\<L\n_ZX,f_*Y\>+\<LX,f_*\n_ZY+\a(Z,Y)\>\\
&\;\;\;+\<B(Z,Y),f_*X\>+\<L\n_ZY,f_*X\>+\<LY,f_*\n_ZX+\a(Z,X)\>\\
&=\<B(Z,X),f_*Y\>+\<LX,\a(Z,Y)\>+\<B(Z,Y),f_*X\>+\<LY,\a(Z,X)\>\\
&=\<(B-f_*\Y\a)(Z,X),f_*Y\>+\<(B-f_*\Y\a)(Z,Y),f_*X\>.
\end{align*}
On the other hand,
$$
\<B(X,Y),f_*Z\>=\<\nab_X\nab_Y\T-\nab_{\nabla_XY}\T,f_*Z\>
-\<X,\nabla\rho\>\<Y,Z\>.
$$
It follows that
$$
C(X,Y,Z)=C(Y,X,Z)\;\;\;\text{and}\;\;\;C(Z,X,Y)=-C(Z,Y,X)
$$
for any $X,Y,Z\in\mathfrak{X}(M)$. Then
\begin{align*}
C(X,Y,Z)&=-C(X,Z,Y)=-C(Z,X,Y)=C(Z,Y,X)\\
&=C(Y,Z,X)=-C(Y,X,Z)=-C(X,Y,Z)=0,
\end{align*}
as we wished.\qed

\begin{proposition}\label{system}
The pair of tensors $(\beta,{\cal E})$ associated to a 
conformal infinitesimal bending  satisfies the system of equations
$\hyperlink{S}{(S)}$.
\end{proposition}

\proof We first show that
\be\label{derY}
(\nab_X\Y)\eta=-f_*B_\eta X-LA_\eta X+{\cal E}(X,\eta)
\ee
where $(\nab_X \Y)\eta=\nab_Xf_*\Y\eta-f_*\Y\nap_X\eta$. 
Taking the derivative of \eqref{Y}, we have from \eqref{inf} and \eqref{Y} that
\begin{align*}
0&=\<\nab_Xf_*\Y\eta,f_*Y\>+\<\Y\eta,\n_XY\>
+\<\nab_XLY,\eta\>+\<LY,\nab_X\eta\>\\
&=\<(\nab_X\Y)\eta,f_*Y\>+\<B_\eta X,Y\>+\<LA_\eta X,f_*Y\>.
\end{align*} 
Since $\<\Y\eta,\xi\>=0$ we have
\begin{align*}
 0&=\<\nab_Xf_*\Y\eta,\xi\>+\<f_*\Y\eta,\nab_X\xi\>
 =\<(\nab_X\Y)\eta,\xi\>-\<\a(X,\Y\eta),\xi\>\\
 &=\<(\nab_X\Y)\eta,\xi\>+\<LA_\eta X-{\cal E}(X,\eta),\xi\>
\end{align*}
for any $X\in\mathfrak{X}(M)$ and $\eta,\xi\in\Gamma(N_fM)$, 
and hence \eqref{derY} follows.

Since
\be\label{form}
(\nab_XB)(Y,Z)
=\nab_X(\nab_YL)Z-(\nab_{\nabla_XY}L)Z-(\nab_Y L)\nabla_XZ
\ee
it is easy to see that 
\be\label{segderL}
(\nab_XB)(Y,Z)-(\nab_YB)(X,Z)=-LR(X,Y)Z
\ee
for all $X,Y,Z\in\mathfrak{X}(M)$. 
It follows using Lemma \ref{tam} that
\begin{align*}
\<(\nab_X B)(Y,Z),f_*W\>&=\<(\nab_X \Y)\a(Y,Z)+f_*\Y(\nap_X \a)(YZ)
-f_*A_{\beta(Y,Z)}X,f_*W\>\\
&\;\;\;+\<Y,W\>\hess\rho(Z,X)-\<Y,Z\>\hess\rho(X,W)
\end{align*}
for any $X,Y,Z,W\in\mathfrak{X}(M)$. Then from
\eqref{segderL}, the Gauss equation 
$$
R(Y,X)Z=A_{\a(X,Z)}Y-A_{\a(Y,Z)}X
$$
and the Codazzi equation, we obtain
\begin{align*}
&\<(\nab_X\Y)\a(Y,Z)-(\nab_Y \Y)\a(X,Z),f_*W\>\\
&=\<LA_{\a(X,Z)}Y-LA_{\a(Y,Z)}X
+f_*A_{\beta(Y,Z)}X-f_*A_{\beta(X,Z)}Y,f_*W\>\\
&+\<Y,Z\>\hess\rho(X,W)-\<Y,W\>\hess\rho(Z,X)
+\<X,W\>\hess\rho(Y,Z)\\&-\<X,Z\>\hess\rho(Y,W).
\end{align*}
On the other hand, it follows from \eqref{derY} that
\begin{align*}
\<(\nab_X \Y)&\a(Y,Z)-(\nab_Y \Y)\a(X,Z),f_*W\>\\
&=\<f_*B_{\a(X,Z)}Y+LA_{\a(X,Z)}Y
-f_*B_{\a(Y,Z)}X-LA_{\a(Y,Z)}X,f_*W\>.
\end{align*}
The last two equations give
\begin{align*}
\<B_{\a(X,Z)}Y&-B_{\a(Y,Z)}X,f_*W\>=
\<A_{\beta(Y,Z)}X-A_{\beta(X,Z)}Y,W\>\\
&+\<Y,Z\>\hess\rho(X,W)-\<Y,W\>\hess\rho(Z,X)\\
&+\<X,W\>\hess\rho(Y,Z)-\<X,Z\>\hess\rho(Y,W),    
\end{align*}
and this is \eqref{derGaussC}.

From Lemma \ref{tam} and \eqref{form} we obtain
\begin{align*}
((\nab_XB)(Y,Z))_{N_fM}&=
\alpha(X,\Y\a(Y,Z))+(\nabla_X^{\perp}\beta)(Y,Z)\\
&\;\;\;+\<Z,\nabla\rho\>\a(X,Y)-\<Y,Z\>\a(X,\nabla\rho).
\end{align*}
Then, we have from \eqref{segderL} and the Gauss equation that
\begin{align*}
(\nabla_X^{\perp}&\beta)(Y,Z)-(\nabla_Y^{\perp}\beta)(X,Z)\\
&=(LR(Y,X)Z)_{N_fM}-\alpha(X,\Y\a(Y,Z))+\alpha(Y,\Y\a(X,Z))\\
&\;\;\;+\<Y,Z\>\a(X,\nabla\rho)-\<X,Z\>\a(Y,\nabla\rho)\\
&=(LA_{\a(X,Z)}Y-LA_{\a(Y,Z)}X)_{N_fM}
-\alpha(X,\Y\a(Y,Z))+\alpha(Y,\Y\a(X,Z))\\
&\;\;\;+\<Y,Z\>\a(X,\nabla\rho)-\<X,Z\>\a(Y,\nabla\rho),
\end{align*}
and this is \eqref{derCodazziC}. 

We have
\begin{align*}
(\nap_X{\cal E})(Y,\eta)&=\nap_X{\cal E}(Y,\eta)
-{\cal E}(\nabla_XY,\eta)-{\cal E}(Y,\nap_X\eta)\\
&=(\nap_X\a)(Y,\Y\eta)+(L(\nabla_X A)(Y,\eta))_{N_fM}
+\a(Y,\nabla_X\Y\eta)\\
&\;\;\;-\a(Y,\Y\nap_X\eta)-(L\nabla_XA_\eta Y)_{N_fM}
+\nap_X(LA_\eta Y)_{N_fM}.
\end{align*}
Then \eqref{derY} yields
\begin{align*}
 (\nap_X{\cal E})(Y,\eta)&=(\nap_X\a)(Y,\Y\eta)
 +(L(\nabla_X A)(Y,\eta))_{N_fM}-\a(Y,B_\eta X)\\
 &\;\;\;-\a(Y,(LA_\eta X)_{TM})-(L\nabla_XA_\eta Y)_{N_fM}
 +\nap_X(LA_\eta Y)_{N_fM}.
\end{align*}

Using the Codazzi equation, we obtain
\begin{align*} 
(\nap_X{\cal E})(Y,\eta)&-(\nap_Y{\cal E})(X,\eta)=
\a(X,B_\eta Y)-\a(Y,B_\eta X)+\a(X,(LA_\eta Y)_{TM})\\
&-\a(Y,(LA_\eta X)_{TM})-(L\nabla_XA_\eta Y)_{N_fM}
+\nap_X(LA_\eta Y)_{N_fM}\\
&+(L\nabla_YA_\eta X)_{N_fM}-\nap_Y(LA_\eta X)_{N_fM}.
\end{align*}
Since
$$
\beta(X, A_\eta Y)=\a(X,(LA_\eta Y)_{TM})-(L\nabla_XA_\eta Y)_{N_fM}
+\nap_X(LA_\eta Y)_{N_fM},
$$
then \eqref{derRicciC} follows.\qed

\section{Trivial infinitesimal bendings}

In this section,  we characterize the trivial conformal 
infinitesimal bendings in terms of the pair of 
tensors $(\beta,{\cal E})$ associated to the bending.
\vspace{1ex}

An infinitesimal bending $\mathcal{T}$ of 
$f\colon M^n\to\R^m$ is trivial if we have $\T=\D f+w$ 
where $\D\in\End(\R^m)$ is skew-symmetric and $w\in\R^m$;
see \cite{DR} or \cite{DT} for details. Then
$L=\D|_{f_*TM}$ and $B(X,Y)=\D\a(X,Y)$. Hence
$$
\beta(X,Y)=\D^N\a(X,Y)\;\;\text{and}\;\;{\cal E}(X,\eta)
=-(\nap_X\D^N)\eta
$$  
where $\D^N\in\Gamma(\End(N_fM))$  given by 
$\D^N\eta=(\D\eta)_{N_fM}$ is skew-symmetric. 

\begin{proposition}\label{trivial}
An infinitesimal bending $\mathcal{T}$ of 
$f\colon M^n\to\R^m$ is trivial if 
and only if there is $C\in\Gamma(\End(N_fM))$ skew-symmetric 
such that 
\be\label{contriv}
\beta(X,Y)=C\a(X,Y)\;\;\text{and}\;\;{\cal E}(X,\eta)
=-(\nap_XC)\eta.
\ee 
\end{proposition}

\proof This is Proposition $5$ in \cite{DJ2}.\vspace{2ex}\qed

It is well-known  that any conformal Killing field on 
an open connected subset of $\R^n$, $n\geq 3$, has the form
$$
X(x)=(\<x,v\>+\lambda)x-(1/2)\|x\|^2v+Cx+w
$$
where  $\lambda\in\R$, $v,w\in\R^n$,  $C\in\End(\R^n)$
is skew-symmetric and the conformal factor is 
$\rho=\<x,v\>+\lambda$; 
cf.\ \cite{Sc} for details.
\medskip

Let $\T$ be a trivial conformal infinitesimal
bending of $f$, that is, locally
$$
\T(x)=(\<f(x),v\>+\lambda)f(x)-1/2\|f(x)\|^2v+\D f(x)+w
$$
where $\lambda\in\R$, $v,w\in\R^m$ and $\D\in\End(\R^m)$ is 
skew-symmetric. Then the conformal factor is 
$\rho(x)=\<f(x),v\>+\lambda$ and 
$$
LX=\<f_*X,v\>f(x)-\<f_*X,f(x)\>v+\D f_*X.
$$
Hence 
\begin{align*}
(\nab_XL)Y=&\,\<f_*Y,v\>f_*X-\<X,Y\>v+\<\a(X,Y),v\>f(x)
-\<\a(X,Y),f(x)\>v\\ &+\D\a(X,Y).
\end{align*}
If $\D'\in \Gamma(\End(f^*T\R^m))$ is the skew-symmetric 
map given by
$$
\D'\sigma=\<\sigma,v\>f(x)-\<\sigma,f(x)\>v+\D \sigma,
$$
then $LX=\D'X$.  Moreover, we have $f_*\Y\eta=(\D'\eta)_{f_*TM}$ and 
$$
(\nab_XL)Y=\<f_*Y,v\>f_*X-\<X,Y\>v+\D'\a(X,Y).
$$
Let $\D^N\in\Gamma(\End(N_fM))$ be given by
$\D^N\xi=(\D'\xi)_{N_fM}$.
Then
$$
\beta(X,Y)=\D^N\a(X,Y)-\<X,Y\>v_N
$$
where $v_N=(v)_{N_fM}$.

We have  
$$
(\nab_X\D')\sigma=\nab_X\D'\sigma-\D'\nab_X\sigma
=\<\sigma,v\>f_*X-\<\sigma,f_*X\>v
$$
for any $X\in\mathfrak{X}(M)$ and $\sigma\in\Gamma(f^*T\R^m)$. 
Then
\begin{align*}
\mathcal{E}(X,\xi)&=\a(X,\Y\xi)+(LA_\xi X)_{N_fM}\\
&=(\nab_X\D'\xi-\nab_{X}\D^N\xi)_{N_fM}+(LA_\xi X)_{N_fM}\\
&=((\nab_X\D')\xi+\D'\nab_X\xi-\nab_X\D^N\xi+LA_\xi X)_{N_fM}\\
&=-(\nap_X\D^N)\xi
\end{align*}
for any $X\in\mathfrak{X}(M)$ and $\xi\in\Gamma(N_fM)$. 

\begin{proposition}\label{trivialC}
A conformal infinitesimal bending $\mathcal{T}$ of 
$f\colon M^n\to\R^m$, \mbox{$n\geq 3$}, is trivial 
if and only if there exist $\delta\in\Gamma(N_fM)$ and 
$C\in\Gamma(\End(N_fM))$ skew-symmetric such that the 
associated  pair has the form
\be\label{contrivC}
\beta(X,Y)=C\a(X,Y)-\<X,Y\>\delta\;\;\text{and}\;\;{\cal E}(X,\xi)
=-(\nap_XC)\xi.
\ee
\end{proposition}

\proof If $(\beta,\mathcal{E})$ has the form \eqref{contrivC} 
and $\rho$ is the conformal factor of $\T$,
we obtain from \eqref{derGaussC} that
\begin{align*}\label{trivgauss}
&\<X,Z\>(\<\a(Y,W),\delta\>-\hess\rho(Y,W))
+\<Y,W\>(\<\a(X,Z),\delta\>-\hess\rho(X,Z))\\
&-\<X,W\>(\<\a(Y,Z),\delta\>\!-\!\hess\rho(Y,Z))
\!-\!\<Y,Z\>(\<\a(X,W),\delta\>\!-\!\hess\rho(X,W))\\
&=0
\end{align*}
for any $X,Y,Z,W\in\mathfrak{X}(M)$.
For $X,Y,W$ orthonormal and $Z=X$ this gives
$$
\<\a(Y,W),\delta\>-\hess\rho(Y,W)=0
$$
whereas for $X=Z$ and $Y=W$ orthonormal this yields
$$
\<\a(X,X),\delta\>-\hess\rho(X,X)
=-\<a(Y,Y),\delta\>+\hess\rho(Y,Y)=0.
$$
Therefore
\be\label{const1}
\<\a(X,Y),\delta\>-\hess\rho(X,Y)=0
\ee
for any $X,Y\in\mathfrak{X}(M)$.

Since $\beta$ and $\mathcal{E}$ have the form \eqref{contrivC} we 
obtain from \eqref{derCodazziC} and the Codazzi equation that
$$
\<X,Z\>(\nap_Y\delta+\a(Y,\nabla\rho))
=\<Y,Z\>(\nap_X\delta+\a(X,\nabla\rho).
$$
Hence
\be\label{const2}
\nap_X\delta+\a(X,\nabla\rho)=0
\ee
for any $X\in\mathfrak{X}(M)$.

Equations \eqref{const1} and \eqref{const2} are equivalent 
to $f_*\nabla\rho+\delta=v$ being constant along $f$.
In particular $\rho(x)=\<f(x),v\>+\lambda$ for 
some $\lambda\in\R$.

Let $\T_1\in\Gamma(f^*T\R^m)$ be the trivial conformal 
infinitesimal bending
$$
\T_1(x)=(\<f(x),v\>+\lambda)f(x)-1/2\|f(x)\|^2v.
$$
Then $\T-\T_1$ is an infinitesimal bending whose 
associated tensors have the form \eqref{contriv}, and thus 
is trivial.\qed

\begin{remark} {\em Two conformal infinitesimal bendings  
$\T_i$, $i=1,2$, of a submanifold  $f\colon M^n\to\R^m$ 
differ by a trivial one if and only if the pairs of 
tensors $(\beta_i,{\cal E}_i)$, $i=1,2$, differ by 
tensors as in \eqref{contrivC}.
}\end{remark}

\begin{corollary}
A conformal infinitesimal bending $\T$ of 
$f\colon M^n\to \R^{n+1}$, $n\geq3$, is trivial if and only 
if its associated tensor ${\cal B}$ has the form
${\cal B}=\varphi I$ for $\varphi\in C^\infty(M)$.
\end{corollary}

\proof In this case, the tensor ${\cal E}$ vanishes and
${\cal B}$ is given by 
$$
\<{\cal B}X,Y\>=\<\beta(X,Y),N\>
$$ 
where $N$ is the Gauss map of $f$. Then $\T$ is trivial 
if and only if 
$$
\beta(X,Y)=-\<X,Y\>\delta
$$
for some $\delta\in\Gamma(N_fM)$. This is equivalent to
${\cal B}=\varphi I$ for $\varphi=-\<\delta,N\>$.
\vspace{2ex}\qed

We conclude this section with some nontrivial examples of 
conformal infinitesimal bendings of rather simple geometric 
nature.

\begin{examples}{\em $(i)$ If $f\colon M^n\to\R^m$ is an  
isometric immersion then a conformal Killing vector field 
of $M^n$ is a conformal infinitesimal bending of~$f$.\vspace{1ex}\\
$(ii)$ Let $g\colon M^n\to\Sf^m$ be an isometric immersion. 
Then ${\cal T}=\varphi f$ is a conformal infinitesimal bending
of $f=i\circ g\colon M^n\to\R^{m+1}$ where $\varphi\in C^\infty(M)$
and $i\colon\Sf^m\to\R^{m+1}$ is the inclusion.
}\end{examples}

\section{The Fundamental theorem}

In this section, we prove the Fundamental theorem for conformal 
infinitesimal variations stated in the introduction.  \vspace{2ex}

Let $\mathbb{V}^{m+1}\subset\Les^{m+2}$ denote the 
\emph{light cone} of the Lorentz space $(\Les^{m+2},\<\,,\,\>)$
defined by
$$
\mathbb{V}^{m+1}=\{v\in\Les^{m+2}\colon\<v,v\>=0,v\neq 0\}.
$$
Given $w\in \mathbb{V}^{m+1}$, then
$$
\mathbb{E}^{m}=\{v\in\mathbb{V}^{m+1}\colon \<v,w\>=1\}
$$
is a model of Euclidean space $\R^m$ in $\Les^{m+2}$.  
In fact, given $v\in\mathbb{E}^{m}$ and a linear isometry 
$C\colon \R^m\to (\spa\{v,w\})^{\perp}\subset\Les^{m+2}$,
the map $\Psi\colon\R^m\to\mathbb{V}^{m+1}\subset\Les^{m+2}$ 
given~by
\be\label{psi}
\Psi(x)=v+Cx-\frac{1}{2}\|x\|^2w
\ee
is an isometric embedding such that 
$\Psi(\R^m)=\mathbb{E}^{m}$.
The normal bundle of $\Psi$ is $N_\Psi\R^m=\spa\{\Psi,w\}$ 
and its second fundamental form is 
\be\label{sff}
\a^{\Psi}(U,V)=-\<U,V\>w
\ee
for any $U,V\in T\R^m$. For further details we refer to 
Section $9.1$ in \cite{DT}.
\vspace{2ex}

\noindent\emph{Proof of Theorem \ref{fundamental}:}
Let $F\colon M^n\to\mathbb{V}^{m+1}\subset\Les^{m+2}$ 
be the isometric immersion $F=\Psi\circ f$, where $\Psi$ 
is given by \eqref{psi}. 
By \eqref{sff} the second fundamental form of $F$  
satisfies
\be\label{sff2}
\a^F(X,Y)=\Psi_*\a(X,Y)-\<X,Y\>w
\ee 
for any $X,Y\in\mathfrak{X}(M)$.

Let $\hat{\beta}\colon TM\times TM\to N_FM$ be 
the symmetric tensor given by
$$
\hat{\beta}(X,Y)=\Psi_*\beta(X,Y)-\hess\rho(X,Y)F
$$
for any $X,Y\in\mathfrak{X}(M)$.
Then \eqref{derGaussC} is equivalent to 
\be\label{derGauss}
A_{\hat{\beta}(Y,Z)}^FX+\hat{B}_{\a^F(Y,Z)}X
-A_{\hat{\beta}(X,Z)}^FY-\hat{B}_{\a^F(X,Z)}Y=0,
\ee
where $A^F_{\xi}$ is the shape operator of $F$ with 
respect to $\xi\in\Gamma(N_FM)$ and $\hat{B}_\xi$ 
is given by
$$
\<\hat{B}_{\xi}X,Y\>=\<\hat{\beta}(X,Y),\xi\>.
$$
Let the tensor 
$\hat{{\cal E}}\colon TM\times N_FM\to N_FM$ be given by
$$
\hat{{\cal E}}(X,\Psi_*\eta)
=\Psi_*{\cal E}(X,\eta)-\<\a(X,\nabla\rho),\eta\>F
\;\;\mbox{for}\;\;\eta\in\Gamma(N_fM),
$$
$$
\hat{{\cal E}}(X,w)=\Psi_*\a(X,\nabla\rho)\;\;\mbox{and}\;\;
\hat{{\cal E}}(X,F)=0.
$$
Since ${\cal E}$ satisfies \eqref{anti} then also does 
$\hat{{\cal E}}$. For simplicity, from now on we just write 
$\eta$ for $\eta\in\Gamma(N_fM)$ as well as its image under 
$\Psi_*$. We have 
$$
(\nabla_X^{'\perp}\hat{\beta})(Y,Z)=(\nap_X\beta)(Y,Z)
-(\nabla_X\hess\rho)(Y,Z)F,
$$
where $\nabla^{'\perp}$ is the normal connection of $F$.
Then
\begin{align*}
(\nabla_X^{'\perp}\hat{\beta})(Y,Z)
-&(\nabla_Y^{'\perp}\hat{\beta})(X,Z)=
(\nap_X\beta)(Y,Z)-(\nap_Y\beta)(X,Z)\\
&+((\nabla_Y\hess\rho)(X,Z)-(\nabla_X\hess\rho)(Y,Z))F.
\end{align*}

It follows from $\hess\rho(X,Y)=\<\nabla_X\nabla\rho,Y\>$
and the Gauss equation that
\begin{align}\label{hess}
(\nabla_Y\hess\rho)(X,Z)
-(\nabla_X&\hess\rho)(Y,Z)=\<R(Y,X)\nabla\rho,Z\>\\
&=\<\a(Y,Z),\a(X,\nabla\rho)\>-\<\a(Y,\nabla\rho),\a(X,Z)\>\nonumber,
\end{align}
where $R$ is the curvature tensor of $M^n$. Thus, 
from \eqref{derCodazziC} and \eqref{hess} we have 
\begin{align*}
(\nabla_X^{'\perp}\hat{\beta})(Y,Z)
&-(\nabla_Y^{'\perp}\hat{\beta})(X,Z)
={\cal E}(Y,\a(X,Z))-{\cal E}(X,\a(Y,Z))\\
&+\<Y,Z\>\a(X,\nabla\rho)-\<X,Z\>\a(Y,\nabla\rho)
+\<\a(Y,Z),\a(X,\nabla\rho)\>F\\
&-\<\a(Y,\nabla\rho),\a(X,Z)\>F,
\end{align*}
and hence
\be\label{derCodazzi}
(\nabla_X^{'\perp}\hat{\beta})(Y,Z)
-(\nabla_Y^{'\perp}\hat{\beta})(X,Z)
=\hat{{\cal E}}(Y,\a^F(X,Z))-\hat{{\cal E}}(X,\a^F(Y,Z))
\ee
for any $X,Y,Z\in\mathfrak{X}(M)$.

From the definition of $\hat{{\cal E}}$ it follows that
$$
(\nabla_X^{'\perp}\hat{{\cal E}})(Y,\eta)
=(\nap_X{\cal E})(Y,\eta)-\<(\nap_X\a)(Y,\nabla\rho),\eta\>F
-\<\a(Y,\nabla_X\nabla\rho),\eta\>F
$$
for any $X,Y\in\mathfrak{X}(M)$ and $\eta\in\Gamma(N_fM)$.
Making use of the Codazzi equation, we obtain
\begin{align*}
(\nabla_X^{'\perp}\hat{{\cal E}})(Y,\eta)
&-(\nabla_Y^{'\perp}\hat{{\cal E}})(X,\eta)
=(\nap_X{\cal E})(Y,\eta)-(\nap_Y{\cal E})(X,\eta)\\
&+(\<\a(X,\nabla_Y\nabla\rho),\eta\>
-\<\a(Y,\nabla_X\nabla\rho),\eta\>)F.
\end{align*}
On the other hand, we have
\begin{align*}
&\hat{\beta}(X,A^F_\eta Y)-\hat{\beta}(A_\eta^FX,Y)
+\a^F(X,\hat{B}_\eta Y)-\a^F(\hat{B}_\eta X,Y)=\beta(X,A_\eta Y)\\
&-\hess\rho(X,A_\eta Y)F-\beta(A_\eta X,Y)+\hess\rho(A_\eta X,Y)F
+\a(X,B_\eta Y)\\
&-\<X,B_\eta Y\>w-\a(B_\eta X,Y)+\<B_\eta X,Y\>w.
\end{align*}
From \eqref{derRicciC}, the symmetry of $\beta$ and 
$$
\hess\rho(X,A_\eta Y)=\<\a(Y,\nabla_X\nabla\rho),\eta\>
$$
it follows that 
\begin{align}\label{derRicci1}
(\nabla_X^{'\perp}\hat{{\cal E}})(Y,\eta)
&-(\nabla_Y^{'\perp}\hat{{\cal E}})(X,\eta)\\
&=\hat{\beta}(X,A^F_\eta X)-\hat{\beta}(A_\eta^FX,Y)
+\a^F(X,\hat{B}_\eta Y)-\a^F(\hat{B}_\eta X,Y)\nonumber
\end{align}

Using the Codazzi equation again, we obtain
$$
(\nabla_X^{'\perp}\hat{{\cal E}})(Y,w)
-(\nabla_Y^{'\perp}\hat{{\cal E}})(X,w)
=\a(Y,\nabla_X\nabla\rho)-\a(X,\nabla_Y\nabla\rho).
$$
Notice that $A^F_w=0$ and $\hat{B}_wX
=-\nabla_X\nabla\rho$ for any $X\in\mathfrak{X}(M)$. 
Then
\begin{align*}
\hat{\beta}(X,A^F_w X)&-\hat{\beta}(A_w^FX,Y)
+\a^F(X,\hat{B}_w Y)-\a^F(\hat{B}_w X,Y)\\
&=\a(Y,\nabla_X\nabla\rho)-\<Y,\nabla_X\nabla\rho\>w
-\a(X,\nabla_Y\nabla\rho)+\<X,\nabla_Y\nabla\rho\>w,
\end{align*}
and hence
\begin{align}\label{derRicci2}
(\nabla_X^{'\perp}\hat{{\cal E}})(Y,w)
&-(\nabla_Y^{'\perp}\hat{{\cal E}})(X,w)\\
&=\hat{\beta}(X,A^F_w Y)-\hat{\beta}(A_w^FX,Y)
+\a^F(X,\hat{B}_w Y)-\a^F(\hat{B}_w X,Y).\nonumber
\end{align}
Since $\hat{B}_F=0$, $A^F_FX=-X$, ${\cal E}(X,F)=0$ 
and $\nabla_X^{'\perp}F=0$, then  
\begin{align}\label{derRicci3}
(\nabla_X^{'\perp}\hat{{\cal E}})(Y,F)
&-(\nabla_Y^{'\perp}\hat{{\cal E}})(X,F)\\
&=\hat{\beta}(X,A^F_FY)-\hat{\beta}(A_F^FX,Y)
+\a^F(X,\hat{B}_FY)-\a^F(\hat{B}_F X,Y)\nonumber
\end{align}
trivially holds.

Summarizing, we have that $\hat{\beta}$ is symmetric, that 
$\hat{{\cal E}}$ satisfies \eqref{anti} and they both verify 
\eqref{derGauss},\eqref{derCodazzi},\eqref{derRicci1}, 
\eqref{derRicci2} and \eqref{derRicci3}.
In this situation, the Fundamental theorem 
for infinitesimal variations, namely, Theorem $6$ 
in \cite{DJ2},  applies.  Notice that in the introduction of
\cite{DJ2} it was observed that this result holds for ambient
spaces of any signature, in particular, for the Lorentzian space
considered here.
Making also use of Proposition $5$ of \cite{DJ2}, we conclude 
that there is an infinitesimal bending 
$\tilde{\T}\in\Gamma(F^*(T\Les^{m+2}))$ of $F$ whose associated 
pair of tensors $(\tilde{\beta},\tilde{\mathcal{E}})$ satisfies
\be\label{tensors}
\tilde{\beta}=\hat{\beta}+\tilde{C}\a^F\;\;\mbox{and}\;\;
\tilde{\mathcal{E}}=\hat{\cal E}-\nabla^{'\perp} \tilde{C}
\ee
where $\tilde{C}\in\Gamma(\End (N_FM))$ is skew-symmetric. 
Moreover, we have that $\tilde{\T}$ is unique up to 
trivial infinitesimal bendings. Write $\tilde{\T}$ as
$$
\tilde{\T}=\Psi_*\T+\<\tilde{\T},w\>F+\<\tilde{\T},F\>w.
$$
Since $\tilde{\T}$ is an infinitesimal bending of $F$, 
we have
$$
\<\nab'_X\tilde{\T},F_*Y\>+\<\nab'_Y\tilde{\T},F_*X\>=0
$$
for all $X,Y\in\mathfrak{X}(M)$. Then
$$
\<\nab_X\T,f_*Y\>+\<\nab_Y\T,f_*X\>+2\<\tilde{\T},w\>\<X,Y\>=0
$$
for all $X,Y\in\mathfrak{X}(M)$.
Hence, setting $\rho_1=-\<\tilde{\T},w\>$ we have that 
$\T$ is a conformal infinitesimal bending of $f$ with 
conformal factor $\rho_1$.

Let $\beta'$ and ${\cal E}'$ be the tensors associated 
to the conformal infinitesimal bending $\T$ of $f$. 
We observe that $(\tilde{\beta})_{\Psi_*N_fM}$ coincides 
with $\beta'$. Let $C\in\Gamma(\End(N_fM))$ be given by 
$C\eta=(\tilde{C}\eta)_{\Psi_*N_fM}$
for any $\eta\in\Gamma(N_fM)$. Then $C$ is skew symmetric. 
It follows from \eqref{sff2} and \eqref{tensors}  that 
the tensor $\beta'$ satisfies
\be\label{tensor1}
\beta'(X,Y)=\beta(X,Y)+C\a(X,Y)-\<X,Y\>\delta
\ee
where $\delta=(\tilde{C}w)_{\Psi_*N_fM}$.

Let $\tilde{L}$ be associated to $\tilde{\T}$ and 
let $\tilde{\Y}$ be given by \eqref{Y} with respect 
to $\tilde{L}$. Given $\eta\in\Gamma(N_fM)$ we have 
$$
\<\tilde{L}X,\eta\>=\<\nab_X\T,\eta\>=\<LX,\eta\>
$$
and then $\tilde{\Y}\eta=\Y\eta$, here $L$ and $\Y$ 
are associated to $\T$. This implies that 
$(\tilde{{\cal E}}(X,\eta))_{\Psi_*N_fM}$ coincides 
with ${\cal E'}(X,\eta)$ for any $X\in\mathfrak{X}(M)$ 
and $\eta\in\Gamma(N_fM)$. Notice that $\Psi_*N_fM$ is 
parallel with respect to $\nabla^{'\perp}$, then we have 
from \eqref{tensors} that
\be\label{tensor2}
{\cal E}'={\cal E}-\nap C.
\ee

Finally it follows from \eqref{tensor1}, \eqref{tensor2} 
and Proposition \ref{trivialC} that any other conformal 
infinitesimal bending arising in this manner differs from 
$\T$ by a trivial conformal infinitesimal bending, and 
this concludes the proof.\qed
\vspace{2ex}

\noindent\emph{Proof of Corollary \ref{corollary}:}
In this case  $\mathcal{E}$ 
vanishes and $\<{\cal B}X,Y\>=\<\beta(X,Y),N\>$. 
Thus \eqref{derRicciC}  holds trivially for $\beta$ 
and $\mathcal{E}=0$. Moreover, 
by the  assumptions on $\mathcal{B}$ we have that $(\beta,0,\rho)$ 
satisfies \eqref{derGaussC} and \eqref{derCodazziC}. Hence 
Theorem \ref{fundamental} gives that $(\beta,0,\rho)$ determines 
a unique conformal infinitesimal bending $\T$ of $f$.\qed

\section{Conformal infinitesimal rigidity}

Let $V^n$ be an $n$-dimensional real vector space and
let $W^{p,p}$ be a real vector space of dimension $2p$
endowed with an indefinite inner product of signature $(p,p)$.
A bilinear form $\gamma\colon V^n\times V^n\to W^{p,p}$ 
is said to be \emph{flat} if
$$
\<\gamma(X,Z),\gamma(Y,W)\>
-\<\gamma(X,W),\gamma(Y,Z)\>=0
$$
for all $X,Y,Z,W\in V^n$.
We say that the bilinear form $\gamma$ is \emph{null} if
$$
\<\gamma(X,Z),\gamma(Y,W)\>=0
$$
for all $X,Y,Z,W\in V^n$. Thus a null bilinear form is flat.
\vspace{1ex}

Given $\gamma\colon V^n\times V^n\to W^{p,p}$ we denote
$$
\mathcal{N}(\gamma)=\{X\in V^n:
\gamma(X,Y)=0\;\;\mbox{for all}\;\;Y\in V^n\}
$$
and
$$
\Sal(\gamma)=\spa\{\gamma(X,Y):X,Y\in V^n\}.
$$

We will need the following result from the theory of
flat bilinear forms which is known to be false for $p=6$;
see Proposition $4.24$ in \cite{DT}.

\begin{lemma}\label{main}
Let $\gamma\colon V^n\times V^n\to W^{p,p}$, $p\leq 5$,
be a symmetric flat bilinear form. 
If $\dim \mathcal{N}(\gamma)\leq n-2p-1$ then there is
an orthogonal decomposition
$$
W^{p,p}=W_1^{\ell,\ell}\oplus W_2^{p-\ell,p-\ell},\;1\leq\ell\leq p,
$$
such that the $W_j$-components $\gamma_j$ of $\gamma$ satisfy:
\begin{itemize}
\item[(i)] $\gamma_1$ is nonzero but null since 
$\Sal(\gamma_1)=\Sal(\gamma)\cap\Sal(\gamma)^\perp$.
\item[(ii)] $\gamma_2$ is flat and
$\dim\mathcal{N}(\gamma_2)\geq n-2p+2\ell$.
\end{itemize}
\end{lemma}

\proof This is Theorem $3$ in \cite{DF} or 
Lemma $4.22$ in \cite{DT}.\qed

\begin{proposition} Let $\T$ be a conformal infinitesimal 
bending of an isometric immersion $f\colon M^n\to\R^m$ with 
conformal factor $\rho$ and associated tensor $\beta$.
Then the bilinear form
$\theta\colon TM\times TM\to N_fM\oplus\R\oplus N_fM\oplus\R$ 
defined at any point of $M^n$ by
{\em \be\label{theta}
\theta=(\alpha+\beta,\<\,,\,\>+\hess\rho,
\alpha-\beta,\<\,,\,\>-\hess\rho)
\ee}
is flat  with respect to the inner product $\<\!\<\,,\,\>\!\>$
of signature $(m-n+1,m-n+1)$ given by
$$
\<\!\<(\xi_1,a_1,\eta_1,b_1),(\xi_2,a_2,\eta_2,b_2)\>\!\>
=\<\xi_1,\xi_2\>_{N_fM}+a_1a_2-\<\eta_1,\eta_2\>_{N_fM}-b_1b_2.
$$
\end{proposition}

\proof A straightforward computation yields
\begin{align*}
&\frac{1}{2}\left(\<\!\<\theta(X,W),\theta(Y,Z)\>\!\>
-\<\!\<\theta(X,Z),\theta(Y,W)\>\!\>\right)
=\<\beta(X,W),\alpha(Y,Z)\>\\
&+\<\alpha(X,W),\beta(Y,Z)\>
-\<\beta(X,Z),\alpha(Y,W)\>-\<\alpha(X,Z),\beta(Y,W)\>\\
&+\<X,W\>\hess\rho(Y,Z)+\<Y,Z\>\hess\rho(X,W)-\<X,Z\>\hess\rho(Y,W)\\
&-\<Y,W\>\hess\rho(X,Z)
\end{align*}
for any $X,Y,Z,W\in\mathfrak{X}(M)$, and the proof follows from 
\eqref{derGaussC}.\qed

\begin{lemma}\label{lema}
Let $f\colon M^n\to\R^m$ be
an isometric immersion.
Let  $Z_1,Z_2\in T_xM$ at $x\in M^n$ be non-zero vectors
such that either $Z_1=Z_2$ or $\<Z_1,Z_2\>=0$.
If $n\geq 4$ and $\nu_1^c(x)\leq n-3$, then
$$
N_fM(x)=\spa\{\a(X,Y)\colon X,Y\in T_xM; 
\<X,Y\>=\<X,Z_1\>=\<Y,Z_2\>=0\}.
$$
\end{lemma}

\proof First assume that $\<Z_1,Z_2\>=0$. 
Let $U^s\subset N_fM(x)$ be the subspace given by
$U^s\perp\a(X,Y)$ for any $X,Y\in T_xM$ as in the statement.
If in addition $\<X,Z_2\>=\<Y,Z_1\>=0$ and $\|X\|=\|Y\|$, 
then $\a(X+Y,X-Y)_{U^s}=0$ gives
$$
\a_{U^s}(X,X)=\a_{U^s}(Y,Y).
$$
Thus, there is $\zeta\in U^s$ such that
$$
\a_{U^s}(X,Y)=\<X,Y\>\zeta
$$
for any $X,Y\in\spa\{Z_1,Z_2\}^\perp$.
Since  $\a_{U^s}(W,Z_1)=\a_{U^s}(W,Z_2)=0$ 
for any $W\in\spa\{Z_1,Z_2\}^\perp$ by assumption, then 
$$
\spa\{Z_1,Z_2\}^\perp\subset\mathcal{N}(\a_{U^s}-\<\,,\,\>\zeta),
$$
and this contradicts our assumption on $\nu_1^c$ unless $s=0$.

If $Z_1=Z_2$ we have again that there is $\zeta\in U^s$ such that
$$
\a_{U^s}(X,Y)=\<X,Y\>\zeta
$$
for any $X,Y\in\spa\{Z_1\}^\perp$. It follows that $A_\zeta$ has an 
eigenspace of multiplicity  at least $n-2$ contradicting the
assumption on  $\nu_1^c$.\qed

\begin{proposition}\label{unique E}
Let $f\colon M^n\to\R^m$, $n\geq 4$, be an isometric immersion and 
let $\T$ be a conformal infinitesimal bending of $f$ with conformal 
factor $\rho$ and associated pair of tensors $(\beta, \mathcal{E})$. 
If $\nu_1^c(x)\leq n-3$ at any $x\in M^n$, then $\mathcal{E}$ is the 
unique tensor satisfying \eqref{anti} as well as an equation of the form
\begin{align}\label{derCodazzigen}
(\nabla^{\perp}_X\beta)(Y,Z)-&(\nabla_Y^{\perp}\beta)(X,Z)
=\mathcal{E}(Y,\a(X,Z))-\mathcal{E}(X,\a(Y,Z))\nonumber\\
&\;\;\;+\<Y,Z\>\psi(X)-\<X,Z\>\psi(Y)
\end{align}
where $\psi\in\Gamma(\Hom(TM,N_fM))$.
\end{proposition}

\proof If also $\mathcal{E}_0\colon TM\times N_fM\to N_fM$ 
satisfies \eqref{anti} and \eqref{derCodazzigen}, then 
\eqref{derCodazziC}  gives
\begin{align*}   
(\mathcal{E}-\mathcal{E}_0)(X,\a(Y,Z))-&(\mathcal{E}
-\mathcal{E}_0)(Y,\a(X,Z))
+\<Y,Z\>(\psi(X)-\a(X,\nabla\rho))\\&
-\<X,Z\>(\psi(Y)-\a(Y,\nabla\rho))=0.
\end{align*}
Then
$$
(\mathcal{E}-\mathcal{E}_0)(X,\a(Y,Z))
=(\mathcal{E}-\mathcal{E}_0)(Y,\a(X,Z))
$$
if $Z$ is orthogonal to $X$ and $Y$. Writing
$$
\<(\mathcal{E}-\mathcal{E}_0)(X_1,\a(X_2,X_3)),\a(X_4,X_5)\>
=(X_1,X_2,X_3,X_4,X_5)
$$
and taking $\<X_1,X_3\>=\<X_2,X_3\>=0$
we have symmetry in the pairs $\{X_1,X_2\}$, $\{X_2,X_3\}$ and
$\{X_4,X_5\}$.  Moreover, since $\mathcal{E}$ and $\mathcal{E}_0$
verify \eqref{anti} we obtain
$$
(X_1,X_2,X_3,X_4,X_5)=-(X_1,X_4,X_5,X_2,X_3).
$$
Hence, if $\{X_i\}_{1\leq i\leq5}$ satisfies 
\be\label{eqs}
\<X_1,X_3\>=\<X_1,X_4\>
=\<X_2,X_3\>=\<X_2,X_5\>=\<X_4,X_5\>=0,
\ee
then
\begin{align*}
(X_1,X_2,X_3,X_4,X_5)&=-(X_1,X_4,X_5,X_2,X_3)=-(X_5,X_4,X_1,X_2,X_3)\\
&=(X_5,X_2,X_3,X_4,X_1)=(X_3,X_2,X_5,X_4,X_1)\\
&=-(X_3,X_4,X_1,X_2,X_5)=-(X_4,X_3,X_1,X_2,X_5)\\
&=(X_4,X_2,X_5,X_3,X_1)=(X_2,X_4,X_5,X_3,X_1)\\
&=-(X_2,X_3,X_1,X_4,X_5)=-(X_1,X_2,X_3,X_4,X_5)\\
&=0.
\end{align*}
Thus
$$
\<(\mathcal{E}-\mathcal{E}_0)(X_1,\a(X_2,X_3)),\a(X_4,X_5)\>=0
$$
if \eqref{eqs} holds. We already have 
$\<X_1,X_4\>=\<X_2,X_5\>=\<X_4,X_5\>=0$. 
Hence, if also $\<X_1,X_2\>=0$ we obtain from Lemma \ref{lema} 
that
$$
(\mathcal{E}-\mathcal{E}_0)(X_1,\a(X_2,X_3))=0
$$
for any $X_1,X_2,X_3\in\mathfrak{X}(M)$ with
$\<X_1,X_2\>=\<X_1,X_3\>=\<X_2,X_3\>=0$.
Then using Lemma \ref{lema} again, it follows that
$$
(\mathcal{E}-\mathcal{E}_0)(X,\eta)=0
$$
for any $X\in \mathfrak{X}(M)$ and $\eta\in\Gamma(N_fM)$.\qed

\begin{lemma}\label{extension}
Let $S\subset\R^m$ be a vector subspace and let 
$T_0\colon S\to \R^m$ be a linear map that is an 
isometry between $S$ and $T_0(S)$. Assume there is no 
$0\neq v\in S$ such that $T_0v=-v$. Then there is 
an isometry $T\in\End(\R^m)$ that extends $T_0$ and 
has $1$ as the only possible real eigenvalue.
\end{lemma}

\proof Extend $T_0$ to an isometry $T$ of $\R^m$. Suppose 
that the eigenspace  of the eigenvalue  $-1$ of $T$
satisfies $\dim E_{-1}=k>0$. We have by assumption that 
$E_{-1}\cap S=E_{-1}\cap T_0(S)=\{0\}$. Let $\{e_1,\ldots,e_k\}$ 
be an orthonormal basis of $E_{-1}$ and set
$$
P=T_0(S)\oplus\spa\{e_2,\ldots,e_k\}.
$$
Let $\xi\in P^\perp$ be a unit vector collinear with the
$P^\perp$-component of $e_1$. Let $\eta\in\R^m$ be such
that $T\eta=\xi$ and let $H$ be the hyperplane $\{\eta\}^\perp$.
If $R$ is the reflection with respect to the hyperplane 
$\{\xi\}^\perp$, then the isometry $T_1=RT$ satisfies
$T_1v=Tv$ for any $v\in H$ since $Tv\in\{\xi\}^\perp$.

Since $\<\eta,e_1\>=-\<\xi,e_1\>\neq 0$, there is 
$v\in H$ such that $\eta+v$ is collinear with $e_1$.  
Hence
$$
T(\eta+v)=\xi+Tv=-\eta-v.
$$
We claim that no vector of the form $\eta+u$, $u\in H$, is 
an eigenvector of $T_1$ associated to $-1$. If otherwise
$$
T_1(\eta+u)=-\xi+Tu=-\eta-u
$$
for some $u\in H$. From the last two equations we obtain
$$
T(u+v)=-2\eta-(u+v).
$$
Then
$$
\|T(u+v)\|^2=4+\|u+v\|^2
$$
which contradicts that $T$ is an isometry and proves the claim.

We have proved that the eigenspace of $T_1$ associated 
to $-1$ is contained in $H$, in fact, it is 
$\spa\{e_2\dots...,e_k\}$. Therefore, by composing $T$ with 
$k$ appropriate reflections we obtain an isometry as in the 
statement.\qed

\begin{proposition}\label{triv}
Let $\T$ be a conformal infinitesimal bending of an isometric 
immersion $f\colon M^n\to\R^m$. If $\T$ is trivial then $\theta$ 
is null. Conversely, if  $\theta$ is null, $n\geq 4$ and 
$\nu_1^c(x)\leq n-3$ at any $x\in M^n$ then $\T$ is trivial.
\end{proposition}

\proof If $\T$ is a trivial conformal infinitesimal bending of 
$f$, then
$$
\T(x)=(\<f(x),v\>+\lambda)f(x)-1/2\|f(x)\|^2v+\D f(x)+w
$$
for some $\lambda\in\R$, $v,w\in\R^m$ and $\D\in\End(\R^m)$ 
skew-symmetric. Since $\rho(x)=\<f(x),v\>+\lambda$, then
$f_*\nabla\rho=v_{TM}$. Thus
\be\label{vconst}
\<\nab_Xv,f_*Y\>=\hess\rho(X,Y)-\<A_{v_{N_fM}}X,Y\>=0
\ee
for any $X,Y\in\mathfrak{X}(M)$. Moreover, we have seen that
$$
\beta(X,Y)=C\a(X,Y)-\<X,Y\>v_{N_fM}
$$
where $C\in\Gamma(\End(N_fM))$ is skew-symmetric.  
Using \eqref{vconst} and that $C$ is skew-symmetric, we obtain
that the bilinear form $\theta$ is null. In fact,
\begin{align*}
\frac{1}{2}\<\!\<\theta(X,Y),\theta(Z,W)\>\!\>
&=\;\<\a(X,Y),\beta(Z,W)\>+\<\beta(X,Y),\a(Z,W)\>\\
&\;\;+\<X,Y\>\hess\rho(Z,W)+\<Z,W\>\hess\rho(X,Y)\\
&=\;-\<Z,W\>\<\a(X,Y),v_{N_fM}\>-\<X,Y\>\<\a(Z,W),v_{N_fM}\>\\
&\;\;+\<X,Y\>\hess\rho(Z,W)+\<Z,W\>\hess\rho(X,Y)\\
&=\,0.
\end{align*}

For the converse, that $\theta$ is null means that
\begin{align}\label{null}
\<\a(X,Y),\beta(Z,W)\>
&+\<\beta(X,Y),\a(Z,W)\>
+\<X,Y\>\hess\rho(Z,W)\nonumber\\&+\<Z,W\>\hess\rho(X,Y)=0
\end{align}
for any $X,Y,Z,W\in\mathfrak{X}(M)$. 
Let $S\subset N_fM(x)\oplus\R$ be the subspace given by
$$
S=\spa\{(\a(X,Y)+\beta(X,Y),\<X,Y\>+\hess\rho(X,Y))\colon X,Y\in T_xM\}.
$$
Then, the map $T_0$ defined by
\begin{align*}
T_0(\a(X,Y)+\beta(X,Y),&\<X,Y\>+\hess\rho(X,Y))\\
&=(\a(X,Y)-\beta(X,Y),\<X,Y\>-\hess\rho(X,Y))
\end{align*}
is an isometry between $S$ and $T(S)$.
We claim that $-1$ is not an eigenvalue of $T_0$. 
If otherwise, then $T_0v=-v$ where 
$$
0\neq v=\sum_i(\a(X_i,Y_i)+\beta(X_i,Y_i),\<X_i,Y_i\>
+\hess\rho(X_i,Y_i))\in S
$$
Hence $\sum_i\a(X_i,Y_i)=0$ and $\sum_i\<X_i,Y_i\>=0$.
Now \eqref{null} gives
$$
\sum_i\<\beta(X_i,Y_i),\a(Z,W)\>+\<Z,W\>\sum_i\hess\rho(X_i,Y_i)=0
$$
for any $Z,W\in\mathfrak{X}(M)$. That is, we have $A_\eta=-hI$ 
where $\eta=\sum_i\beta(X_i,Y_i)$ and $h=\sum_i\hess\rho(X_i,Y_i)$.
From our assumption on $\nu_1^c$ we obtain $\eta=h=0$, hence $v=0$ 
proving the claim.

Let $T$ be the isometry of $N_fM(x)\oplus\R$ extending
$T_0$ given by Lemma~\ref{extension}. Then
$$
(I+T^t)(I-T)=(T^t-T)=-(I-T^t)(I+T)
$$
where $T^t$ is the transpose of $T$. Thus
$$
(I-T)(I+T)^{-1}=-(I+T^t)^{-1}(I-T^t)=-((I-T)(I+T)^{-1})^t, 
$$
that is,  $(I-T)(I+T)^{-1}$ is a skew-symmetric endomorphism 
of $N_fM(x)\oplus\R$.

It is easy to see that
$$
(I-T)(I+T)^{-1}(\a(X,Y),0)=(\beta(X,Y),\hess\rho(X,Y))
$$
for any $X,Y\in T_xM$ such that $\<X,Y\>=0$. 
Thus, there is $C\in\End(N_fM(x))$ skew-symmetric such that 
$\beta(X,Y)=C\a(X,Y)$ for any $X,Y\in T_xM$ with $\<X,Y\>=0$.
It follows from $\beta(X+Y,X-Y)=C\a(X+Y,X-Y)$ for any $X,Y\in T_xM$ 
orthonormal that
$$
\beta(X,X)-C\a(X,X)=\beta(Y,Y)-C\a(Y,Y).
$$
Therefore, there is $\delta\in N_fM(x)$ such that 
\be\label{betatriv}
\beta(X,Y)=C\a(X,Y)-\<X,Y\>\delta
\ee
for any $X,Y\in T_xM$.

By Lemma \ref{lema} there are smooth local vector fields 
$X_i,Y_i$, $1\leq i\leq m-n$, satisfying $\<X_i,Y_i\>=0$  
such that $\a(X_i,Y_i)$, $1\leq i\leq m-n$, span the normal 
bundle. Thus $C$ and $\delta$ are smooth. 

Define $\mathcal{E}_0\colon TM\times N_fM\to N_fM$ by 
$\mathcal{E}_0(X,\eta)=-(\nap_X C)\eta$. It follows from \eqref{betatriv} that
\begin{align*}
    (\nap_X\beta)(Y,Z)-(\nap_Y\beta)(X,Z)&=\mathcal{E}_0(Y,\a(X,Z))-\mathcal{E}_0(X,\a(Y,Z))\\
&\;\;-\<Y,Z\>\nap_X\delta+\<X,Z\>\nap_Y\delta.
\end{align*}
Then we have from Lemma \ref{unique E} that $\mathcal{E}=\mathcal{E}_0$, 
and thus $\T$ is trivial by Proposition \ref{trivialC}.
\vspace{2ex}\qed

\noindent\emph{Proof of Theorem \ref{rigidity}:} 
Let $\T$ be a conformal infinitesimal bending of $f$
such that the flat bilinear $\theta$ given by \eqref{theta}
is not null at $x\in M^n$.  Since 
$\mathcal{N}(\theta)=\{0\}$, there is an orthogonal decomposition 
$$
W_0^{2p+2}=N_fM(x)\oplus\R\oplus N_fM(x)\oplus\R
=W_1^{\ell,\ell}\oplus W_2^{p-\ell+1,p-\ell+1},\; 1\leq\ell\leq p,
$$
such that $\theta$ splits as $\theta=\theta_1+\theta_2$
as in Lemma \ref{main}. Denoting
$\Delta=\mathcal{N}(\theta_2)$, we have $\dim\Delta\geq n-2(p-\ell+1)$.
Thus $\theta(Z,X)=\theta_1(Z,X)$
for any $Z\in\Delta$ and $X\in T_xM$.

Let $S\subset N_fM(x)\oplus\R$ be the vector subspace given by
$$
S=\spa\{(\a(Z,X)+\beta(Z,X),\<Z,X\>
+\hess\rho(Z,X))\colon Z\in\Delta\,\mbox{and}\,X\in T_xM\}.
$$
If $\Pi_1$ denotes the orthogonal projection from $W_0^{2p+2}$ 
onto the first copy of $N_fM(x)\oplus\R$, 
then $S\subset\Pi_1(\Sal(\theta)\cap\Sal(\theta)^\perp)$
and, in particular, $\dim S\leq\ell$.

That $\theta_1$ is null means that the map 
$T\colon S\to N_fM(x)\oplus\R$ defined by
\begin{align*}
T(\a(Z,X)+&\beta(Z,X),\<Z,X\>+\hess\rho(Z,X))\\
&=(\a(Z,X)-\beta(Z,X),\<Z,X\>-\hess\rho(Z,X)).
\end{align*}
is an isometry between $S$ and $T(S)$. We have that
$$
\frac{1}{2}(I+T)(\a(Z,X)+\beta(Z,X),\<Z,X\>
+\hess\rho(Z,X))=(\a(Z,X),\<Z,X\>).
$$
If $S_1=((I+T)(S))^\perp\subset N_fM\times\R$, 
then $\dim S_1\geq p-\ell+1$.
For $(\eta,a)\in S_1$
\be\label{i}
\<\a(Z,X),\eta\>+a\<X,Z\>=0
\ee
for any $Z\in\Delta$ and $X\in T_xM$.
Let $U\subset N_fM$ be the orthogonal projection of $S_1$ 
in $N_fM$. Since $S_1$ does not posses elements of the form 
$(0,a)$ with $0\neq a\in\R$, then
$$
\dim U\geq p-\ell+1.
$$
It follows from \eqref{i} that there  exists $\zeta\in U$ such that
$$
\a_U(Z,X)=\<Z,X\>\zeta
$$
for any $Z\in\Delta$ and $X\in T_xM$. Hence $\a_U-\<\,,\,\>\zeta$ 
has a kernel of dimension  at least $\dim\Delta\geq n-2(p-\ell+1)$. 
But this contradicts the assumption on the conformal $s$-nullities, 
and hence $\theta$ is necessarily null at any point. 
We conclude from Proposition \ref{triv} that $\T$ is trivial.\qed

\section*{Acknowledgments}

This work is the result of the visit 21171/IV/19 funded by the
Fundación Séneca-Agencia de Ciencia y Tecnología de la
Región de Murcia in connection with the ``Jiménez De La Espada"
Regional Programme For Mobility, Collaboration And Knowledge Exchange.

Marcos Dajczer was partially supported by the Fundación Séneca Grant\\
21171/IV/19 (Programa Jiménez de la Espada), MICINN/FEDER project
PGC2018-097046-B-I00, and Fundación Séneca project 19901/GERM/15, 
Spain.

Miguel I. Jimenez thanks the mathematics department of the Universidad 
de Murcia for the hospitality during his visit where part of this 
work was developed.

\noindent Marcos Dajczer\\
IMPA -- Estrada Dona Castorina, 110\\
22460--320, Rio de Janeiro -- Brazil\\
e-mail: marcos@impa.br

\bigskip

\noindent Miguel Ibieta Jimenez\\
Instituto de Ciências Matemáticas e de Computação\\
Universidade de São Paulo\\
São Carlos\\
SP 13560-970-- Brazil\\
e-mail: mibieta@impa.br

\begin{thebibliography}{lll}

\bibitem{Ca0} Cartan, E., 
\emph{La d\'eformation des hypersurfaces dans l'espace euclidien 
r\'eel a $n$ dimensions},
Bull. Soc. Math. France {\bf 44} (1916), 65--99.

\bibitem{Ca}  Cartan, E., \emph{La d\'eformation des
hypersurfaces dans l'espace conforme r\'eel a $n\geq 5$ dimensions},
Bull. Soc. Math. France {\bf 45} (1917), 57--121.

\bibitem{Ce} Cesàro, E., ``Lezioni di Geometria intrinseca''.
Napoli, 1896.

\bibitem{CD} do Carmo, M. and Dajczer, M., 
\emph{Conformal rigidity}, 
Amer. J. Math. {\bf 109} (1987),  963--985. 

\bibitem{DF} Dajczer, M. and Florit, L., 
\emph{Compositions of isometric immersions in higher
codimension},
Manuscripta Math. {\bf 105} (2001), 507--517. 

\bibitem{DJ} Dajczer, M. and Jimenez, M. I., 
\emph{Genuine infinitesimal bendings of submanifolds}, 
preprint. 


\bibitem{DJ2} Dajczer, M. and Jimenez, M. I., 
\emph{Infinitesimal variations of submanifolds}, 
to appear in Bull. Braz. Math. Soc. 


\bibitem{DJV} Dajczer M., Jimenez M. I. and Vlachos, Th.,
\emph{Conformal infinitesimal variations of Euclidean 
hypersurfaces}, preprint.


\bibitem{DR} Dajczer, M. and Rodríguez, L., 
\emph{Infinitesimal rigidity of Euclidean submanifolds}, 
Ann. Inst. Fourier {\bf 40} (1990), 939--949.

\bibitem{DT0} Dajczer, M. and Tojeiro, R.,
\emph{On Cartan's conformally deformable hypersurfaces},
Michigan Math. J. {\bf 47} (2000), 529--557.

\bibitem{DT} Dajczer, M. and  Tojeiro, R., 
``Submanifold theory beyond an introduction".
Universitext. Springer, 2019.

\bibitem{DV} Dajczer, M. and Vlachos, Th.,
\emph{The infinitesimally bendable Euclidean hypersurfaces},
Ann. Mat. Pura Appl. {\bf 196} (2017), 1961--1979 and 
Ann. Mat. Pura Appl. {\bf 196} (2017), 1981--1982.

\bibitem{Ji} Jimenez, M. I., 
\emph{Infinitesimal bendings of complete Euclidean 
hypersurfaces}, 
Manuscripta Math. {\bf 157} (2018), 513--527.

\bibitem{Sb} Sbrana, U.,
\emph{Sulla deformazione infinitesima delle ipersuperficie},
Ann. Mat. Pura  Appl. {\bf 15} (1908), 329--348.

\bibitem{Sb2} Sbrana, U.,
\emph{Sulle variet\`a ad $n-1$ dimensioni deformabili nello 
spazio euclideo ad $n$ dimensioni},
Rend. Circ. Mat. Palermo  {\bf 27} (1909), 1--45.

\bibitem{Sc} Schottenloher, M., 
``A mathematical introduction to conformal field theory".
Lecture Notes in Physics, 759. Springer-Verlag, 2008.

\bibitem{To} Tojeiro, R.,
\emph{Liouville's theorem revisited},
Enseign. Math. {\bf 53} (2007), 67--86.

\bibitem{Ya} Yano, K. 
\emph{Notes on infinitesimal variations of submanifolds}, 
J. Math. Soc. Japan {\bf 32} (1980), 45--53. 

\end{thebibliography}
\end{document}